\newcommand{\R}{\mathds R}

\newcommand{\Z}{\mathds Z}

\newcommand{\Ker}{\mathrm{Ker}}

\newcommand{\Bsym}{\mathrm{B_{\mathrm{sym}}}}
\newcommand{\indice}{\mathrm n_-}
\newcommand{\coindice}{\mathrm n_+}
\newcommand{\segnatura}{\mathrm{sign}}
\newcommand{\Dim}{\mathrm{dim}}

\newcommand{\Gr}{\mathrm{Gr}}

\newcommand{\Spl}{\mathrm{Sp}}

\newcommand{\spl}{\mathrm{sp}}

\newcommand{\iCZ}{\mathfrak i_{\mathrm{CZ}}}

\newcommand{\Lsa}{\mathcal L_{\mathrm{sa}}}

\documentclass[oneside,draft,a4paper,10pt]{amsart}

\usepackage{times}
\usepackage{dsfont}

\title[A product formula for the Conley--Zehnder index]{On a product formula for the Conley--Zehnder index of
symplectic paths and its applications}
\author[M. de Gosson]{Maurice de Gosson}
\author[S. de Gosson]{Serge de Gosson}
\author[P.\ Piccione]{Paolo Piccione}
\address{Departamento de Matem\'atica\\ Instituto de Matem\'atica e Estat\'\i stica\\ Universidade de S\~ao Paulo\\
Rua do Mat\~ao 1010, CEP 05508-090, S\~ao Paulo, SP, Brazil}

\date{July 1st, 2006}

\begin{document}


\theoremstyle{plain}\newtheorem{teo}{Theorem}[section]
\theoremstyle{plain}\newtheorem{prop}[teo]{Proposition}
\theoremstyle{plain}\newtheorem{lem}[teo]{Lemma}
\theoremstyle{plain}\newtheorem{cor}[teo]{Corollary}
\theoremstyle{definition}\newtheorem{defin}[teo]{Definition}
\theoremstyle{remark}\newtheorem{rem}[teo]{Remark}
\theoremstyle{plain} \newtheorem{assum}[teo]{Assumption}
\swapnumbers
\theoremstyle{definition}\newtheorem{example}{Example}[section]
\theoremstyle{plain} \newtheorem*{acknowledgement}{Acknowledgements}


\begin{abstract}
Using invariance by fixed-endpoints homotopies and a generalized notion
of symplectic Cayley transform, we prove a product formula for the Conley--Zehnder index of
continuous paths  with arbitrary endpoints in the symplectic group. We discuss two applications of the formula,
to the metaplectic group and to periodic solutions of Hamiltonian systems.
\end{abstract}

\maketitle

\begin{section}{Introduction}
The theory of periodic Hamiltonian orbits plays a fundamental role in many
active parts of both pure and applied mathematics. An object of choice for the
qualitative and quantitative study of these orbits is the Conley--Zehnder
index, introduced in \cite{ConZeh}, and whose theory has been further studied
by many authors, see e.g. \cite{CapLeeMil,HWZ,RobSal} and the references
therein. Let us explain briefly what this index is about. Let $H\in C^{\infty
}(\mathbb{R}^{2n}\times\mathbb{R},\mathbb{R})$ be a time-dependent
Hamiltonian, and\ denote by $(\mathcal F_{t}^{H})$ the flow it determines; assume that
$z_{0}\in\mathds{R}^{2n}$ is such that $\mathcal F_{1}^{H}(z)=z$, the mapping
$t\longmapsto \mathcal F_{t}^{H}(z)$ is a $1$-periodic orbit through $z_{0}%
\in\mathds{R}^{2n}$. Assume now that the Jacobian matrix $\Phi(z)=D\mathcal F_1(z)$ of
$\mathcal F_{1}^{H}$ satisfies the non-degeneracy condition%
\begin{equation}
\det\big(\Phi(z)-\mathrm{Id}\big)\neq0;\label{2}%
\end{equation}
Set now $z(t)=\mathcal F_{t}^{H}(z)$ and consider the linearized Hamiltonian system
along $z$; its time-evolution is governed by the ordinary
differential equation $\dot{u}=JD^{2}H(z,t)u$ whose flow consists of the
symplectic matrices $\Phi_{t}=D\mathcal F_{t}^{H}(z)$. The path $\Sigma:t\longmapsto
\Phi_{t}$, $t\in\left]0,1\right]$ lies in the symplectic group $\Spl(\R^{2n},\omega)$;
it starts from the identity and ends at $\Phi(z)$. If the
non-degeneracy condition (\ref{2}) holds one associates to $\Sigma$ an integer
$\iCZ(\Sigma)$, the \emph{Conley-Zehnder index} of the path $\Sigma$,
and whose vocation is to give an algebraic count of the number of points
$t_{j}$ in the interval $\left]0,1\right[$ for which $\Phi(z)-\mathrm{Id}$ is not invertible.  One should
however be aware of the fact that condition (\ref{2}) is very restrictive; in
particular it is never satisfied in the simple case where $H$ is
time-independent! The aim of this paper is to give a general definition of the
Conley--Zehnder index and to prove a formula for the index of the product of
two symplectic paths, obtaining as a consequence a formula for the Conley--Zehnder
index of an iterated periodic orbit. We observe that Cushman and Duistermaat \cite{cushdu}
and Duistermaat~\cite{Duistermaat} also have addressed the question of the
index of the iteration of periodic orbits; the methods these authors use are
however very different from ours. An extensive literature on the Maslov index
and its iteration properties has been produced by Y. Long and his collaborators
(see for instance \cite{LongAdv}), who obtained remarkable results on the multiplicity of periodic orbits
of Hamiltonian systems.

In order to obtain the product formula in the general case of arbitrary endpoints,
in this paper we introduce the notion of generalized symplectic Cayley transform of a symplectomorphism.
For each symplectomorphism $\psi$ whose spectrum does not contain $1$, we define
a real-analytic diffeomorphism $\mathcal C_\psi$ between (an open subset of) the symplectic
group and (an open subset of) the space of symmetric operators on a real
finite dimensional Hilbert space. The classical symplectic Cayley transform is
obtained when $\psi=-\mathrm{Id}$. The generalized Cayley transform is used to compute
the correction term in the product formula (formula \eqref{eq:prodformula})
of symplectic paths.
\medskip

This article is structured as follows:
We begin by giving a working definition of the usual Conley--Zehnder
index; we take the opportunity to recall a few basic definitions and results
about some well-known objects from symplectic geometry such as the
H\"{o}rmander, Wall--Kashiwara, Maslov, and Leray indices; a good reference is
the seminal paper \cite{CapLeeMil}. The last part of this Section is devoted
to the introduction and study of the main properties of a notion of symplectic
Cayley transform which will be instrumental to our study of the product formula.

 In Section~\ref{sec:prodformula} we state and prove the product formula
for the Conley--Zehnder index in a very general setting in terms of the
H\"{o}rmander index; in the special case where both paths are non-degenerate
this formula can be restated very simply using the symplectic Cayley transform
previously defined. The basic argument employed in the proof of the product formula
uses a  homotopy properties of paths in topological groups

In Section~\ref{sec:iterationformulas} we pursue our study of the
product formula and discuss the problem of the calculation of the index an
orbit which is iterated an arbitrary number of times; we obtain a number of
precise estimates. These results aim at applications in the theory of periodic
Hamiltonian orbits and at applications to spectral flow formulae and Morse theory
(see \cite{ConZeh}, \cite{RobSal2}, \cite{SalZeh}).

 Finally, in Section~\ref{sec:metaplectic} we apply our results to the
Weyl representation of metaplectic operators; this question is of a
fundamental importance in the study of the semiclassical quantization of
non-integrable Hamiltonian systems (Gutzwiller's theory, see \cite{Gutz,WM,muratore}); in
particular we improve previous results \cite{mdglettmath}.

\begin{acknowledgement}
The first author (MdG) has been financed by a FAPESP grant, and he wishes to
thank all faculties and staff at the University of S\~ao Paulo for providing excellent working
conditions during his stay.
\end{acknowledgement}
\end{section}

\begin{section}{Preliminaries on Maslov, Conley--Zehnder, Kashiwara, Leray  and H\"ormander's
indexes}
Let $V$ be a finite dimensional real vector space. By $\Bsym(V)$ we mean the space
of all symmetric bilinear forms $B:V\times V\to\R$  on $V$;
we will always identify a bilinear form $B:V\times V\to\R$
with the linear operator $V\ni v\mapsto B(v,\cdot)\in V^*$, that
will be denoted by the same symbol $B$.
For $B\in\Bsym(V)$, we denote by $\indice(B)$,
$\coindice(B)$ respectively the {\em index} and  the {\em
coindex\/}  of $B$.
The \emph{signature}  of $B$ is the
difference $\segnatura(B)=\indice(B)-\coindice(B)$.
A bilinear form $B$ will be called \emph{nondegenerate} if the
linear map $V\ni v\mapsto B(v,\cdot)\in V^*$ is an isomorphism.

A \emph{symplectic form} on $V$ is a nondegenerate antisymmetric bilinear form
$\omega:V\times V\to\R$; the standard example of a \emph{symplectic space}
is $V=\R^n\oplus{\R^n}^*$ endowed with the \emph{canonical symplectic form}:
\begin{equation}\label{eq:omega0}
\omega_0\big((v,\alpha),(w,\beta)\big)=\beta(v)-\alpha(w),\qquad v,w\in\R^n,\ \alpha,\beta\in{\R^n}^*.
\end{equation}
The {\em symplectic group\/}
$\Spl(V,\omega)$ is the closed subgroup of $\mathrm{GL}(V)$
consisting of those linear maps on $V$ that preserve $\omega$.

Given a symplectic space $(V,\omega)$, with $\Dim(V)=2n$,
a {\em Lagrangian subspace\/} of $V$ is an
$n$-dimensional subspace $L\subset V$ on which
$\omega$ vanishes. The set of all Lagrangian
subspaces of $V$, denoted by $\Lambda=\Lambda(V,\omega)$, has the structure
of a compact, real-analytic submanifold of the Grassmannian
of all $n$-dimensional subspaces of $V$. The dimension
of $\Lambda$ equals $\frac12n(n+1)$, and a real-analytic
atlas on $\Lambda$ is given as follows.

For all $L\in\Lambda$ and $k\in\{0,\ldots,n\}$, let
$\Lambda_0(L)$ denote the set of all Lagrangian
subspaces that are transverse to $L$, which is a dense
open subset of $\Lambda$. Given a pair $L_0,L_1\in\Lambda$
of complementary Lagrangians, i.e., $L_0\cap L_1=\{0\}$, then
one defines a map:
\[\varphi_{L_0,L_1}:\Lambda_0(L_1)\longrightarrow\Bsym(L_0)\]
as follows. Any Lagrangian $L\in\Lambda_0(L_1)$ is the graph
of a unique linear map $T:L_0\to L_1$; then, $\varphi_{L_0,L_1}$
is defined to be the restriction of the bilinear map
$\omega(T\cdot,\cdot)$ to $L_0\times L_0$. It is easy to see
that, due to the fact that $L$ is Lagrangian, such bilinear map is symmetric.

Given $L_1\in\Lambda$ and $L_0,L\in\Lambda_0(L_1)$, the
bilinear forms $\varphi_{L_0,L_1}(L)\in\Bsym(L_0)$ and
$\varphi_{L,L_1}(L_0)\in\Bsym(L)$ are related by the identity:
\[\eta^*\varphi_{L,L_1}(L_0)=-\varphi_{L_0,L_1}(L),\]
where $\eta:L_0\to L$ is the isomorphism given by the restriction to
$L_0$ of the projection $\pi^L:L\oplus L_1\cong V\to L$.
In particular:
\begin{equation}\label{eq:minussignature}\segnatura\big(\varphi_{L,L_1}(L_0)\big)=-\segnatura\big(\varphi_{L_0,L_1}(L)\big).\end{equation}
We need another identity relating the charts on $\Lambda$; assume that $L_0,L_1,L,L'$ are
Lagrangian subspaces of $V$, with $L_0,L,L'$ transversal to $L_1$. Then:
\[\varphi_{L,L_1}(L')=\eta_0^*\big(\varphi_{L_0,L_1}(L')-\varphi_{L_0,L_1}(L)\big),\]
where $\eta_0:L\to L_0$ is the isomorphism given by the restriction to $L$ of the projection
$L_0\oplus L_1\cong V\to L_0$. Hence:
\begin{equation}\label{eq:sommasignatures}
\segnatura\big(\varphi_{L,L_1}(L')\big)=\segnatura\big(\varphi_{L_0,L_1}(L')-\varphi_{L_0,L_1}(L)\big).
\end{equation}

Given Lagrangians $L_0,L_1\subset V$ with
$L_0\cap L_1=\{0\}$, then there exists a symplectic isomorphism
(symplectomorphism)
$\phi:V\to\R^n\oplus{\R^n}^*$ (i.e., the pull-back $\phi^*\omega_0$
coincides with $\omega$)
such that
$\phi(L_0)=\{0\}\oplus {\R^n}^*$ and
$\phi(L_1)=\R^n\oplus\{0\}$.
\subsection{Maslov index,  Wall--Kashiwara's index and Leray index}
Denote by $\pi(\Lambda)$ the \emph{fundamental groupoid} of $\Lambda$, i.e., the set of fixed-endpoints homotopy
classes $[\gamma]$ of continuous paths $\gamma$ in $\Lambda$, endowed with the partial operation
of concatenation $\diamond$.
For all $L_0\in\Lambda$, there exists a unique $\frac12\Z$-valued
groupoid homomorphism
$\mu_{L_0}$ on $\pi(\Lambda)$ such that:
\begin{equation}\label{eq:muL0tilde}
\mu_{L_0}\big([\gamma]\big)=\tfrac12\segnatura\big(\varphi_{L_0,L_1}(\gamma(1))\big)-
\tfrac12\segnatura\big(\varphi_{L_0,L_1}(\gamma(0))\big)
\end{equation}
for all continuous curve $\gamma:[0,1]\to\Lambda_0(L_1)$ and for all
$L_1\in\Lambda_0(L_0)$.
\begin{defin}
The map $\mu_{L_0}:\pi(\Lambda)\to\frac12\Z$ is called the \emph{$L_0$-Maslov index}.
\end{defin}
The Maslov index has the following property \emph{(symplectic
invariance)}:
 if $\phi:(V,\omega)\to(V',\omega')$ is a symplectomorphism and
$\gamma:[a,b]\to\Lambda(V,\omega)$ is continuous, then:
\[\mu_{L_0}(\gamma)=\mu_{\phi(L_0)}\big(\phi\circ\gamma\big).\]
Moreover, replacing the symplectic form $\omega$ by $-\omega$ produces
a change in the sign  of $\mu_{L_0}$.

Let $(L_{0},L_{1},L_{2})$ be a triple of elements of $\Lambda$; the
\emph{Wall--Kashiwara index} (see \cite{CapLeeMil}) of that triple is the signature
$\tau(L_{0},L_{1},L_{2})$ of the quadratic form \[L_{0}\oplus L_{1}\oplus L_{2}\ni (z_{0},z_{1},z_{2}
)\longmapsto\omega(z_{0},z_{1})+\omega(z_{1},z_{2})+\omega(z_{2},z_{0})\in\R.\]
It is a $\Spl(V,\omega)$-invariant totally
antisymmetric $2$-cocycle on $\Lambda$. Let $\pi
:\Lambda_{\infty}\longrightarrow\Lambda$ be the universal covering of
$\Lambda$; the \emph{Leray index} on $\Lambda_{\infty}$ is the unique $\mathds Z$-valued
$1$-cochain $\mu$ on $\Lambda_{\infty}$ which is locally constant on
$\{(L_{\infty},L_{\infty}^{\prime}):\pi(L_{\infty})\cap\pi(L_{\infty}^{\prime
})=0\}$ and such that $\partial\mu=\pi^{\ast}\tau$ ($\partial$ the \v{C}ech
coboundary operator). It is a symplectic invariant, in the sense that:
\[
\mu(\Phi_{\infty}L_{\infty},\Phi_{\infty}L_{\infty}^{\prime})=\mu(L_{\infty
},L_{\infty}^{\prime})
\]
for all $\Phi_{\infty}\in\Spl_{\infty}(V,\omega)$ (the universal
covering group of $\Spl(V,\omega)$). The Leray and Maslov
indices are related in the following way: identifying $L_{\infty}\in
\Lambda_{\infty}$ with the fixed-endpoints homotopy classes of continuous
paths $\gamma$ joining $L_{0}$ to $L$, we have (see \cite{CM})
\begin{equation}
\mu_{L_{0}}\big([\gamma]\big)=\tfrac{1}{2}\big[\mu(\gamma(1)_{\infty},L_{0,\infty}%
)-\mu(\gamma(0)_{\infty},L_{0,\infty})\big]\label{gosson1}%
\end{equation}
where $L_{0,\infty}$ is the homotopy class of any loop through $L_{0}$,
$\gamma(0)_{\infty}$ the homotopy class of any path $\gamma_{0}$ joining
$L_{0}$ to $\gamma(0)$, and $\gamma(1)_{\infty}$ that of the concatenation
$\gamma_{0}\diamond\gamma$. In view of the cochain relation $\partial\mu
=\pi^{\ast}\tau$ this formula can be rewritten%
\begin{equation}
\mu_{L_{0}}\big([\gamma]\big)=\tfrac{1}{2}\big[\mu\big(\gamma(1)_{\infty},\gamma(0)_{\infty
}\big)+\tau\big(L_{0},\gamma(0),\gamma(1)\big)\big].\label{gosson2}%
\end{equation}
(See Cappell \textit{et al.} \cite{CapLeeMil} for a comparative
study of Leray and related indexes).
\subsection{H\"ormander's index}
Given four Lagrangians $L_0,L_1,L_0',L_1'\in\Lambda$ and
any continuous curve $\gamma:[a,b]\to\Lambda$ such that
$\gamma(a)=L_0'$ and $\gamma(b)=L_1'$, then the value of
the quantity $\mu_{L_1}(\gamma)-\mu_{L_0}(\gamma)$ does
{\em not\/} depend on the choice of $\gamma$.
\begin{defin}\label{thm:defHormander}
Given  $L_0,L_1,L_0',L_1'\in\Lambda$, the
{\em H\"ormander index\/} $\mathfrak q(L_0,L_1;L_0',L_1')$ is
the half-integer number $\mu_{L_1}(\gamma)-\mu_{L_0}(\gamma)$, where
$\gamma:[a,b]\to\Lambda$ is any continuous curve with
$\gamma(a)=L_0'$ and $\gamma(b)=L_1'$.
\end{defin}

It follows from \eqref{gosson2}, the property $\partial \mu =\pi ^{\ast
}\tau $, and the fact that $\mu (L_{\infty },L_{\infty }^{\prime })$ is
independent of the choice of base point in $\Lambda $ that we have%
\begin{multline*}
-\mathfrak q\big(L_0,L_1;\gamma(0),\gamma(1)\big)=\mu _{L_{0}}\big([\gamma ]\big)-\mu _{L_{1}}\big([\gamma ]\big)\\=
\tfrac{1}{2}\big[\tau
\big(L_{0},\gamma (0),\gamma (1)\big)-\tau \big(L_{1},\gamma (0),\gamma (1)\big)\big].
\end{multline*}
\subsection{The Conley--Zehnder index}
Given a symplectic space $(V,\omega)$, consider the direct
sum $V^2=V\oplus V$, endowed with the symplectic form $\omega^2=\omega\oplus(-\omega)$,
defined by:
\[\omega^2\big((v_1,v_2),(w_1,w_2)\big)=\omega(v_1,v_2)-\omega(w_1,w_2),\quad v_1,v_2,w_1,w_2\in V.\]
Given a linear operator $T:V\to V$, we will denote by $\Gr(T)\subset V^2$ its graph.
Let $\Delta\subset V^2$ denote the diagonal;
if $\Phi\in\Spl(V,\omega)$, then $\Gr(\Phi)=(\mathrm{Id}\oplus\Phi)[\Delta]\in\Lambda(V^2,\omega^2)$;
in particular $\Delta=\Gr(\mathrm{Id})$ and $\Delta^o=\{(v,-v):v\in V\}=\Gr(-\mathrm{Id})$ are Lagrangian subspaces
of $V^2$.
If $\Phi_1,\Phi_2\in\Spl(V,\omega)$, then $\Phi_1\oplus\Phi_2:V^2\to V^2$
belongs to $\Spl(V^2,\omega^2)$.
\begin{defin}\label{thm:defmindsimplgr}
Given a continuous curve $\Phi$ in
$\Spl(V,\omega)$, the \emph{Conley--Zehnder index $\iCZ(\Phi)$ of $\Phi$}
 is the $\Delta$-Maslov
index of the curve $t\mapsto\Gr\big(\Phi(t)\big)\in\Lambda(V^2,\omega^2)$:
\[\iCZ(\Phi):=\mu_\Delta\big(t\mapsto\Gr(\Phi(t))\big).\]
\end{defin}

The above is one of the possible definitions of the notion of Conley--Zehnder index
(see \cite{ConZeh, RobSal, SalZeh}). The Conley--Zehnder index is additive by concatenation
and invariant by fixed endpoint homotopies.

Consider the map $\mathcal A:V^2\to V^2$ given by $\mathcal A(v_1,v_2)=(v_2,v_1)$;
$\mathcal A$ is an \emph{anti-symplectomorphism} of $(V^2,\omega^2)$, i.e., $\mathcal A^*\omega^2=-\omega^2$.
Clearly, $\mathcal A=\mathcal A^{-1}$ and $\mathcal A[\Delta]=\Delta$; more generally, if $S:V\to V$ is
a bijection, then $\mathcal A[\Gr(S)]=\Gr(S^{-1})$.
It follows that, given a continuous path $\Phi:[a,b]\to\Spl(V,\omega)$, one
has:
\begin{multline*}\iCZ(\Phi^{-1})=\mu_\Delta\big(t\mapsto\Gr(\Phi^{-1})\big)=\mu_\Delta\big(t\mapsto\mathcal A\big[\Gr(\Phi(t))\big]\big)\\=
-\mu_{\mathcal A[\Delta]}\big(t\mapsto\Gr(\Phi(t))\big)=-\mu_\Delta\big(t\mapsto\Gr(\Phi(t))\big)=-\iCZ(\Phi).
\end{multline*}
Moreover, given four Lagrangians $L_0,L_1,L_0',L_1'\in\Lambda$:
\[\mathfrak q(\mathcal A[L_0],\mathcal A[L_1];\mathcal A[L_0'],\mathcal A[L_1'])=-\mathfrak q(L_0,L_1;L_0',L_1').\]
\begin{lem}\label{thm:iCZprodq}
Let $\Phi:[a,b]\to\Spl(V,\omega)$ be a continuous curve, and let
$\psi_*\in\Spl(V,\omega)$ be fixed. Denote by $\psi_*\cdot\Phi$ and by $\Phi\cdot\psi_*$ the
continuous curves in $\Spl(V,\omega)$ given by $t\mapsto\psi_*\cdot\Phi(t)$ and $t\mapsto\Phi(t)\cdot\psi_*$
respectively. Then:
\[
\iCZ(\psi_*\cdot\Phi)=\iCZ(\Phi\cdot\psi_*)=\iCZ(\Phi)+\mathfrak q\Big(\Delta,\Gr(\psi_*^{-1});\Gr\big(\Phi(a)\big),
\Gr\big(\Phi(b)\big)\Big).\]
\end{lem}
\begin{proof}
A direct computation, as follows:
\begin{multline*}
\iCZ(\psi_*\cdot\Phi)=\mu_\Delta\big(t\mapsto\Gr(\psi_*\Phi(t))\big)=\mu_\Delta\big(t\mapsto(\mathrm{Id}\oplus\psi_*)
(\mathrm{Id}\oplus\Phi(t)))[\Delta]\big)\\
\stackrel{\text{symplectic invariance}}=\mu_{(\mathrm{Id}\oplus\psi_*^{-1})[\Delta]}\big(t\mapsto(\mathrm{Id}\oplus\Phi(t)))[\Delta]\big)
\\=\mu_\Delta\big(t\mapsto(\mathrm{Id}\oplus\Phi(t)))[\Delta]\big)+\mathfrak q\big(\Delta,(\mathrm{Id}\oplus\psi_*^{-1})[\Delta];
\Gr\big(\Phi(a)\big),\Gr\big(\Phi(b)\big)\big)\\=\iCZ(\Phi)+\mathfrak q\Big(\Delta,\Gr(\psi_*^{-1});\Gr\big(\Phi(a)\big),
\Gr\big(\Phi(b)\big)\Big).
\end{multline*}
Similarly,
\begin{multline*}
\iCZ(\Phi\cdot\psi_*)=-\iCZ(\psi_*^{-1}\cdot\Phi^{-1})\\=-\Big[\iCZ(\Phi^{-1})+\mathfrak q\Big(\Delta,\Gr(\psi_*);\Gr\big(\Phi(a)^{-1}\big),\Gr\big(\Phi(b)^{-1}\big)\Big)\Big]
\\=\iCZ(\Phi)+\mathfrak q\Big(\mathcal A[\Delta],\mathcal A[\Gr(\psi_*)];\mathcal A\big[\Gr\big(\Phi(a)^{-1}\big)\big],\mathcal A\big[\Gr\big(\Phi(b)^{-1}\big)\big]\Big)
\\=\iCZ(\Phi)+\mathfrak q\Big(\Delta,\Gr(\psi_*^{-1});\Gr\big(\Phi(a)\big),
\Gr\big(\Phi(b)\big)\Big).\qedhere
\end{multline*}
\end{proof}

\subsection{The symplectic Cayley transform}
Let us consider the two transverse Lagrangians $\Delta,\Delta^o\in\Lambda(V^2,\omega^2)$.
Assume that $\Phi\in\Spl(V,\omega)$ does not have the eigenvalue $1$, i.e.,
that $\mathrm{Id}-\Phi$ is invertible or, equivalently, that $\Gr(\Phi)$ is
transverse to $\Delta$.

An immediate computation shows that
$\varphi_{\Delta^o,\Delta}\big(\Gr(\Phi)\big)$ is identified with
the symmetric bilinear form $2\omega\big((\mathrm{Id}+\Phi)(\mathrm{Id}-\Phi)^{-1}\cdot,\cdot)$
on the vector space $V$. In particular, if $J$ is a complex structure on $V$ and
$\langle\cdot,\cdot\rangle$ is a positive definite inner product on $V$ with $\langle J\cdot,\cdot\rangle=\omega$,
then the linear operator $J(\mathrm{Id}+\Phi)(\mathrm{Id}-\Phi)^{-1}:V\to V$
is symmetric relatively to $\langle\cdot,\cdot\rangle$.
In what follows, we will assume that the vector space $V$ is endowed with a positive definite inner product
$\langle\cdot,\cdot\rangle$ and with a complex structure $J$ as above.
We will denote by $\Lsa(V)$ the space of symmetric linear operators on $V$ and we will implicitly identify the spaces
$\Bsym(V)$ and $\Lsa(V)$ by the obvious identification.
\begin{defin}\label{thm:defCaytransform}
Given objects $V$, $\omega$, $J$ and $\langle\cdot,\cdot\rangle$ as
above, and given $\Phi\in\Spl(V,\omega)$ with $(\mathrm{Id}-\Phi)$ invertible,
then the \emph{symplectic Cayley transform} $M_\Phi$ of $\Phi$ is the
symmetric operator $\frac12J(\mathrm{Id}+\Phi)(\mathrm{Id}-\Phi)^{-1}$.
\end{defin}
One checks immediately the equality: $M_{\Phi^{-1}}=-M_\Phi$.
\smallskip

The notion of symplectic Cayley transform was originally introduced
by  Mehlig and Wilkinson (see \cite{WM}) and further studied in \cite{mdglettmath}
and in \cite{PhD}. We observe that the sign convention used in this paper differs
from the original one.

In order to deal with symplectomorphisms whose spectrum contains $1$,
we need to introduce a generalization of the notion of Cayley transform.

For a fixed $\psi\in\Spl(V,\omega)$, we denote by $\Spl_\psi(V,\omega)$ the
dense open subset of $\Spl(V,\omega)$ consisting of those $\Phi\in\Spl(V,\omega)$
such that $\Phi-\psi$ is invertible. When $\psi=\mathrm{Id}$, it is customary
to write $\Spl_{\mathrm{Id}}(V,\omega)=\Spl_0(V,\omega)$.
\begin{lem}\label{thm:Caleytrasnfwelldefined}
Let $\psi\in\Spl_0(V,\omega)$ be fixed. For all $\Phi\in\Spl_\psi(V,\omega)$, the linear
operator
\begin{equation}\label{eq:defpsiCayley}
\mathcal C_\psi(\Phi)=J(\psi-\mathrm{Id})(\Phi-\psi)^{-1}(\Phi-\mathrm{Id}):V\longrightarrow V
\end{equation}
is symmetric. Moreover, $\Ker\big(\mathcal C_\psi(\Phi)\big)=\Ker(\Phi-\mathrm{Id})$.
\end{lem}
\begin{proof}
The assumption $\psi\in\Spl_0(V,\omega)$ says that $\Gr(\psi)\in\Lambda(V^2,\omega^2)$ is
transverse to $\Delta$; the assumption that $\Phi\in\Spl_\psi(V,\omega)$ says that $\Gr(\Phi)\in\Lambda(V^2,\omega^2)$ is
transverse to $\Gr(\psi)$.
A direct calculation shows that the symmetric bilinear form \[\varphi_{\Delta,\Gr(\psi)}\big(\Gr(\Phi)\big)\in\Bsym(\Delta)\]
can be identified (via the isomorphism $V\ni v\mapsto (v,v)\in\Delta$) with the
bilinear form $\omega\big((\mathrm{Id}-\psi)(\Phi-\psi)^{-1}(\mathrm{Id}-\Phi)\cdot,\cdot\big)=\langle\mathcal C_\psi(\Phi)\cdot,\cdot\rangle$
on $V$. From this observation, the conclusion follows.
The last statement in the thesis is obvious.
\end{proof}

\begin{defin}
Let $\psi\in\Spl_0(V,\omega)$ be fixed. The map $\mathcal C_\psi:\Spl_\psi(V,\omega)\to\Lsa(V)$
defined in Lemma~\ref{thm:Caleytrasnfwelldefined} is called the \emph{$\psi$-Cayley transform}.
\end{defin}

\end{section}

\begin{section}{The product formula}\label{sec:prodformula}
Let us start with a simple result on signature of the difference of symmetric bilinear forms
(a similar result is proven in \cite{RobSal}):
\begin{lem}\label{thm:segnsimmforms}
Let $V$ be a finite dimensional real vector space and let $U,Z\in\Bsym(V)$
be nondegenerate symmetric bilinear forms on $V$ such that $U-Z$ is also
nondegenerate. Then,  $U^{-1}-Z^{-1}$ is nondegenerate and:
\[
\segnatura(Z)-\segnatura(U)=\segnatura(Z^{-1}-U^{-1})-\segnatura(U-Z).
\]
\end{lem}
\begin{proof}
Define the nondegenerate symmetric
bilinear form $B\in\Bsym(V^2)$ by \[B((a_1,b_1),(a_2,b_2))=
Z(a_1,a_2)-U(b_1,b_2).\] Identifying $V$ with $\Delta$ by $v\mapsto(v,v)$, one computes
easily $B\vert_\Delta=Z-U$, which is nondegenerate.
Denote by $\Delta^{\perp_B}$ the $B$-orthogonal complement of $\Delta$;
identifying $V$ with $\Delta^{\perp_B}$ by $V\ni V\to
(v,U^{-1}Zv)\in\Delta^{\perp_B}$, it is easily seen that
$B\vert_{\Delta^{\perp_B}}=Z(Z^{-1}-U^{-1})Z$. The conclusion
follows.
\end{proof}

\begin{prop}\label{thm:calcHormindex}
Let $(V,\omega)$ be a symplectic vector space, and let $L,L',L_0,L_1\in\Lambda(V,\omega)$
be four Lagrangians, with $L,L',L_1$ transverse to $L_0$. Then:
\[\mathfrak q(L_0,L;L_0,L')=\tfrac12\,\segnatura\big(\varphi_{L_1,L_0}(L)-\varphi_{L_1,L_0}(L')\big).\]
\end{prop}
\begin{proof}
Up to a symplectic isomorphism, we can assume $V=\R^{n}\oplus{\R^n}^*$, with $\omega$ the canonical symplectic
form $\omega_0$ (see \eqref{eq:omega0}), $L_0=\{0\}\oplus{\R^n}^*$ and $L_1=\R^n\oplus\{0\}$;
transversality of $L$ and $L'$ with $L_0$ says that we can write $L=\Gr(T)$, $L'=\Gr(T')$, where
$T,T':\R^n\to{\R^n}^*$ are self-adjoint linear maps. By definition, the H\"ormander's index
$\mathfrak q(L_0,L;L_0,L')$ is given by:
\[\mathfrak q(L_0,L;L_0,L')=\mu_L(\gamma)-\mu_{L_0}(\gamma),\]
where $\gamma:[a,b]\to\Lambda$ is an arbitrary continuous curve with $\gamma(a)=L_0$ and $\gamma(b)=L'$.

In order to compute the two Maslov indexes in the formula above, let us choose a
Lagrangian $\widetilde L_1$ which is transverse simultaneously to the three
Lagrangians $L_0$, $L$ and $L'$. We can choose, for instance, $\widetilde L_1=\Gr(S)$,
where $S:\R^n\to{\R^n}^*$ is a self-adjoint linear operator; in this way, $\widetilde L_1$
is transverse to $L_0$. Transversality of $\widetilde L_1$ to $L$ and $L'$ is equivalent to
$T-S$ and $T'-S$ being invertible.
Since $L_0,L'\in\Lambda_0(\widetilde L_1)$, and $\Lambda_0(\widetilde L_1)$ is arc-connected, then
one can choose a curve $\gamma$ with the required properties whose image
is contained in $\Lambda_0(\widetilde L_1)$. Using \eqref{eq:muL0tilde}, we then get:
\begin{multline*}\mathfrak q(L_0,L;L_0,L')\\=\tfrac12\big[\segnatura\big(\varphi_{L,\widetilde L_1}(L')\big)-
\segnatura\big(\varphi_{L,\widetilde L_1}(L_0)\big)-\segnatura\big(\varphi_{L_0,\widetilde L_1}(L')\big)
+\segnatura\big(\varphi_{L_0,\widetilde L_1}(L_0)\big)\big].
\end{multline*}
A direct calculation gives:
\begin{itemize}
\item $\varphi_{L_1,L_0}(L)\cong T:\R^n\to{\R^n}^*$;
\item $\varphi_{L_1,L_0}(L')\cong T':\R^n\to{\R^n}^*$;
\item $\varphi_{L_0,\widetilde L_1}(L_0)=0$;
\item $\varphi_{L_0,\widetilde L_1}(L')\cong ({T'}-S)^{-1}:{\R^n}^*\to\R^n$;
\item $\varphi_{L,\widetilde L_1}(L_0)$ is identified\footnote{%
Here, $L=\Gr(T)$ is identified with $\R^n$ via the map: $\R^n\ni v\mapsto(v,Tv)\in L$.}
with the symmetric bilinear form
on $\R^n$ given by the self-adjoint linear operator $S-T:\R^n\to{\R^n}^*$;
\item $\varphi_{L,\widetilde L_1}(L')$ is identified with the symmetric bilinear form
on $\R^n$ given by:  \[(S-T)+(T-S)(S-T')^{-1}(S-T):\R^n\to{\R^n}^*.\]
\end{itemize}
Using Lemma~\ref{thm:segnsimmforms} and
keeping in mind that $\segnatura\big((S-T)+(T-S)(S-T')^{-1}(S-T)\big)=\segnatura\big((S-T)^{-1}-(S-T')^{-1}\big)$,
from the equalities above one obtains:
\[\mathfrak q(L_0,L;L_0,L')=\tfrac12\segnatura(T-T'),\]
which concludes the proof.
\end{proof}
\begin{lem}\label{thm:prodhomconc} Given continuous paths $\Phi_i:[0,1]\to\Spl(V,\omega)$, $i=1,2$,
then the pointwise product $\Phi(t)=\Phi_1(t)\cdot\Phi_2(t)$ is fixed-endpoints
homotopic to the concatenation $\widetilde\Phi\diamond\widetilde\Phi_2$, where $\widetilde\Phi_i:[0,1]\to\Spl(V,\omega)$
is given by:
\[\widetilde\Phi_1(t)=\Phi_1(t)\cdot\Phi_2(0),\qquad\widetilde\Phi_2(t)=\Phi_1(1)\cdot\Phi_2(t),\qquad\forall\,t\in[0,1].\]
In particular, $\iCZ(\Phi)=\iCZ(\widetilde\Phi_1)+\iCZ(\widetilde\Phi_2)$.
\end{lem}
\begin{proof}
The curve $\Phi_1$ is fixed-endpoints homotopic to the curve $\overline\Phi_1$ defined
by: \[\overline\Phi_1(t)=\begin{cases}\Phi_1(2t),& \text{if $t\in[0,\tfrac12]$};\\
\Phi_1(1),& \text{if $t\in\left[\tfrac12,1\right]$},\end{cases}\]
while $\Phi_2$ is fixed-endpoints homotopic to $\overline\Phi_2$, given by:
\[\overline\Phi_2(t)=\begin{cases}\Phi_2(0),& \text{if $t\in[0,\tfrac12]$};\\
\Phi_2(2t-1),& \text{if $t\in\left[\tfrac12,1\right]$},\end{cases}\]
hence the pointwise product $\Phi_1\cdot\Phi_2$ is fixed-endpoints homotopic to the
pointwise product $\overline\Phi_1\cdot\overline\Phi_2$. Clearly, $\overline\Phi_1\cdot\overline\Phi_2=
\widetilde\Phi_1\diamond\widetilde\Phi_2$.
\end{proof}

An immediate application of Lemma~\ref{thm:iCZprodq} and Lemma~\ref{thm:prodhomconc} gives:

\begin{cor}\label{thm:homtriloop}
Given continuous paths $\Phi_i:[0,1]\to\Spl(V,\omega)$, $i=1,2$,
then:
\begin{multline}\label{eq:firstprodform}\iCZ(\Phi_1\cdot\Phi_2)=\iCZ(\Phi_1)+\iCZ(\Phi_2)\\+\mathfrak q\big(\Delta,\Gr\big(\Phi_2(0)^{-1}\big);\Gr\big(\Phi_1(0)\big),
\Gr\big(\Phi_1(1)\big)\big)\\+\mathfrak q\big(\Delta,\Gr\big(\Phi_1(1)^{-1}\big);\Gr\big(\Phi_2(0)\big),
\Gr\big(\Phi_2(1)\big)\big).\end{multline}
In particular, if $\Phi_2$ is a homotopically trivial loop starting at the identity of $\Spl(V,\omega)$,
then $\iCZ(\Phi_1\cdot\Phi_2)=\iCZ(\Phi_1)$.
\end{cor}

Let us now consider the case that the paths $\Phi_1$ and $\Phi_2$ start at the identity of $\Spl(V,\omega)$, in which case
obviously $\Gr\big(\Phi_1(0)\big)=\Gr\big(\Phi_2(0)\big)=\Delta$, and the term \[\mathfrak q\big(\Delta,\Gr\big(\Phi_2(0)^{-1}\big);\Gr\big(\Phi_1(0)\big),
\Gr\big(\Phi_1(1)\big)\big)=\mathfrak q\big(\Delta,\Delta;\Delta,
\Gr\big(\Phi_1(1)\big)\big)\]  in equality \eqref{eq:firstprodform} vanishes.
\begin{cor}\label{thm:prodformula}
Let $\Phi_i:[0,1]\to\Spl(V,\omega)$ be continuous paths with $\Phi_1(0)=\Phi_2(0)=\mathrm{Id}$ and
with $\mathrm{Id}-\Phi_1(1)$, $\mathrm{Id}-\Phi_2(1)$ invertible. Then:
\begin{equation}\label{eq:prodformula}
\iCZ(\Phi_1\cdot\Phi_2)=\iCZ(\Phi_2\cdot\Phi_1)=\iCZ(\Phi_1)+\iCZ(\Phi_2)-\tfrac12\segnatura\big(M_{\Phi_1(1)}+M_{\Phi_2(1)}\big).
\end{equation}
\end{cor}
\begin{proof}
Formula \eqref{eq:prodformula} follows easily from \eqref{eq:firstprodform}, using Proposition~\ref{thm:calcHormindex}
applied to the symplectic space $(V^2,\omega^2)$ and to the Lagrangians $L_0=\Delta$, $L_1=\Delta^o$, $L=\Gr\big(\Phi_1(1)^{-1}\big)$ and $L'=\Gr\big(\Phi_2(1)\big)$.
\end{proof}

The result of Corollary~\ref{thm:prodformula} can be extended to the case of paths
with arbitrary endpoints in the symplectic group, using our generalized notion of symplectic Cayley transform.
\begin{prop}\label{thm:calcHorgraph}
Let $\phi_1,\phi_2,\psi\in\Spl(V,\omega)$ be fixed; assume that $\psi\in\Spl_0(V,\omega)$
and that $\phi_1,\phi_2\in\Spl_\psi(V,\omega)$. Then, the H\"ormander index $\mathfrak q\big(\Delta,\Gr(\phi_1);\Delta,\Gr(\phi_2)\big)$
is given by:
\[\mathfrak q\big(\Delta,\Gr(\phi_1);\Delta,\Gr(\phi_2)\big)=\tfrac12\Big[\segnatura\big(\mathcal C_\psi(\phi_2)-\mathcal C_\psi(\phi_1)\big)
-\segnatura\big(\mathcal C_\psi(\phi_2)\big)+\segnatura\big(\mathcal C_\psi(\phi_1)\big)\Big].\]
If $\phi_1,\phi_2\in\Spl_0(V,\omega)$, then:
\[\mathfrak q\big(\Delta,\Gr(\phi_1);\Delta,\Gr(\phi_2)\big)=\tfrac12\segnatura\big(\mathcal C_\psi(\phi_1)^{-1}-\mathcal C_\psi(\phi_2)^{-1}\big).\]
\end{prop}
\begin{proof}
By definition of H\"ormander index, $\mathfrak q\big(\Delta,\Gr(\phi_1);\Delta,\Gr(\phi_2)\big)$ is given by:
\[\mu_{\Gr(\phi_1)}(\gamma)-\mu_\Delta(\gamma),\]
where $\gamma:[a,b]\to\Lambda(V^2,\omega^2)$ is any continuous curve with $\gamma(a)=\Delta$ and $\gamma(b)=\Gr(\phi_2)$.
We can choose one such curve $\gamma$ whose image remains inside the set of Lagrangians of $(V^2,\omega^2)$ that
are transversal to $\Gr(\psi)$. Then, by definition of Maslov index, we have:
\begin{multline*}
\mu_{\Gr(\phi_1)}(\gamma)-\mu_\Delta(\gamma)=\tfrac12\Big[\segnatura\Big(\varphi_{\Gr(\phi_1),\Gr(\psi)}\big(\Gr(\phi_2)\big)\Big)-\segnatura\big(
\varphi_{\Gr(\phi_1),\Gr(\psi)}(\Delta)\big)\Big]\\
-\tfrac12\Big[\segnatura\Big(\varphi_{\Delta,\Gr(\psi)}\big(\Gr(\phi_2)\big)\Big)-\segnatura\big(\varphi_{\Delta,\Gr(\psi)}(\Delta)\big)\Big].\end{multline*}
The first equality in the thesis is now obtained easily, using the following:
\smallskip

\begin{itemize}
\item $\varphi_{\Delta,\Gr(\psi)}\big(\Gr(\phi_2)\big)=\mathcal C_\psi(\phi_2)$;
\smallskip

\item $\segnatura\big(\varphi_{\Gr(\phi_1),\Gr(\psi)}(\Delta)\big)\stackrel{\text{by \eqref{eq:minussignature}}}{=}-\segnatura\big(\varphi_{\Delta,\Gr(\psi)}\big(\Gr(\phi_1)\big)
\big)=-\segnatura\big(\mathcal C_\psi(\phi_1)\big)$;
\smallskip

\item $\segnatura\!\Big(\!\varphi_{\Gr(\phi_1),\Gr(\psi)}\big(\Gr(\phi_2)\big)\!\Big)\!\stackrel{\text{by \eqref{eq:sommasignatures}}}=\!\segnatura\!\Big(\!
\varphi_{\Delta,\Gr(\psi)}\big(\Gr(\phi_2)\big)-\varphi_{\Delta,\Gr(\psi)}\big(\Gr(\phi_1)\big)\!\Big)$\hfill\break
\phantom{a}\hfill $=\segnatura\big(\mathcal C_\psi(\phi_2)-\mathcal C_\psi(\phi_1)\big)$;
\smallskip

\item $\varphi_{\Delta,\Gr(\psi)}(\Delta)=0$.
\end{itemize}
\smallskip

Finally, if $\phi_1,\phi_2\in\Spl_0(V,\omega)$, i.e., if $\mathcal C_\psi(\phi_1)$ and $\mathcal C_\psi(\phi_2)$ are invertible,
then, by Lemma~\ref{thm:segnsimmforms}:
\[\segnatura\big(\mathcal C_\psi(\phi_2)-\mathcal C_\psi(\phi_1)\big)
-\segnatura\big(\mathcal C_\psi(\phi_2)\big)+\segnatura\big(\mathcal C_\psi(\phi_1)\big)=\segnatura\big(\mathcal C_\psi(\phi_1)^{-1}-\mathcal C_\psi(\phi_2)^{-1}\big),\]
which concludes the proof.
\end{proof}

\begin{cor}\label{thm:prodformula2}
Let $\Phi_i:[0,1]\to\Spl(V,\omega)$ be continuous paths with $\Phi_1(0)=\Phi_2(0)=\mathrm{Id}$, and let
$\psi\in\Spl_0(V,\omega)$ be such that $\Phi_1(1),\Phi_2(1)\in\Spl_\psi(V,\omega)$. Then:
\begin{multline}\label{eq:prodformula2}
\iCZ(\Phi_1\cdot\Phi_2)=\iCZ(\Phi_1)+\iCZ(\Phi_2)+\tfrac12\segnatura\left[\mathcal C_\psi\big(\Phi_2(1)\big)-\mathcal C_\psi\big(\Phi_1(1)\big)\right]\\
-\tfrac12\segnatura\big[\mathcal C_\psi\big(\Phi_2(1)\big))\big]+\tfrac12\segnatura\big[\mathcal C_\psi\big(\Phi_1(1)\big))\big].
\end{multline}
If $\Phi:[0,1]\to\Spl(V,\omega)$ is a continuous path, then:
\begin{equation}\label{eq:estimateiCZPhiNbis}
\Big\vert \vert\iCZ(\Phi^N)\vert-N\vert\iCZ(\Phi)\vert\Big\vert\le\tfrac12n(N-1),
\end{equation}
for all $N\ge1$.
\end{cor}
\begin{proof}
The proof of formula \eqref{eq:prodformula2} is analogous to the proof of \eqref{eq:prodformula},
where one replaces the chart $\varphi_{\Delta^o,\Delta}$ by $\varphi_{\Delta,\Gr(\psi)}$, using the
result of Proposition~\ref{thm:calcHorgraph} for the computation of the H\"ormander index.
The inequality in formula \eqref{eq:estimateiCZPhiNbis} is obtained readily from \eqref{eq:prodformula2} using induction on $N$;
for such induction argument one has to keep in mind the following observations.
\begin{itemize}
\item One has to use Cayley transform $\mathcal C_\psi$, where $\psi\in\Spl_0(V,\omega)$ is such that
all the powers $\Phi(1)^N$ belong to $\Spl_\psi(V,\omega)$. The set of such $\psi$'s is non empty, and in fact \emph{dense}
in $\Spl(V,\omega)$; namely, this set is the intersection of the countable family of dense open subsets $\Spl_0(V,\omega)\cap\Spl_{\Phi(1)^j}(V,\omega)$,
$j\ge1$, and the claims follows from Baire's theorem.
\item Given any pair $B_1,B_2$ of symmetric bilinear form on any $n$-dimensional real vector space, then
\[\vert\segnatura(B_2-B_1)-\segnatura(B_2)-\segnatura(B_1)\vert\le n.\]
Namely, if $B_1$ and $B_2$
are nondegenerate, the claim follows immediately from Lemma~\ref{thm:segnsimmforms}; for the general
case simply use an argument of density and continuity.
\end{itemize}
\end{proof}
The inequality in formula \eqref{eq:estimateiCZPhiNbis} tells us that if $\big\vert\iCZ(\Phi)\big\vert\ge\frac n2$, then
$\big\vert\iCZ(\Phi^{N})\vert$ has a linear growth in $N$.

\end{section}

\begin{section}{Iteration formulas}\label{sec:iterationformulas}

Let us now discuss the problem of determining the Conley--Zehnder index
and the Maslov index of the iteration of a periodic solution of a
Hamiltonian system.

Let $(\mathcal M,\varpi)$ be a $2n$-dimensional symplectic manifold, and let $H:\mathcal M\times\R\to\R$
be a time-dependent smooth Hamiltonian. Assume that $H$ is $T$-periodic in time,
and that $z:[0,T]\to M$ is a solution of $H$ (i.e., $\dot z=\vec H(z)$
such that $z(0)=z(T)$,
where $\vec H$ is the time-dependent Hamiltonian vector field, defined by $\varpi(\vec H,\cdot)=\mathrm dH$).
Then, the iterates $z^{(N)}$ of $z$, defined as
the concatenation:
\[z^{(N)}=\underbrace{z\diamond\cdots\diamond z}_{\text{$N$-times}}:[0,NT]\longrightarrow\mathcal M\]
are also solutions of $H$. Assume that it is given a \emph{periodic
symplectic trivialization} of the tangent bundle of $\mathcal M$
along $z$ (i.e., of the pull-back $z^*T\mathcal M$), which consists
of a smooth family $\Psi=\{\psi_t\}_{t\in[0,T]}$ of
symlectomorphisms $\psi_t:T_{z(0)}\mathcal M\to T_{z(t)}\mathcal M$
with $\psi_0=\psi_T=\mathrm{Id}$. By a simple orientability
argument, periodic symplectic trivializations along periodic
solutions always exist. By the periodicity assumption, we have a
smooth extension $\R\ni t\mapsto\psi_t$ by setting
$\psi_{t+NT}=\psi_t$ for all $t\in[0,T]$.

Denote by $\mathcal F^H_{t,t'}:\mathcal M\to\mathcal M$ the (maximal) flow of $\vec H$,\footnote{%
For our purposes, we will not be interested in questions of global existence of the flow
$\mathcal F^H$.}
i.e., $\mathcal F^H_{t,t'}(p)=\gamma(t')$, where $\gamma$ is the unique
integral curve of the time-dependent vector field $\vec H$ on $\mathcal M$
satisfying $\gamma(t)=p$. It is well known that for all $t,t'$, the $\mathcal F^H_{t,t'}$ is
a symplectomorphism between open subsets of $\mathcal M$. Conjugation with $\psi_t$ gives
a smooth map $\R\ni t\mapsto X(t)=\psi_t^{-1}\circ\mathcal F^H_{0,t}\big(z(0)\big)\circ\psi_t$
of linear endomorphisms of $T_{z(0)}\mathcal M$; clearly $X(t)$ lies in the Lie algebra
$\spl\big(T_{z(0)}\mathcal M,\varpi_{z(0)}\big)$ of the symplectic group $\Spl\big(T_{z(0)}\mathcal M,\varpi_{z(0)}\big)$.

The \emph{linearized Hamilton equation along $z$} is the linear system
\begin{equation}\label{eq:linhamequations}
v'(t)=X(t)v(t),
\end{equation}
in $T_{z(0)}\mathcal M$; the fundamental solution of this linear system
is a smooth symplectic path $\Phi:\R\to\Spl\big(T_{z(0)}\mathcal M,\varpi_{z(0)}\big)$ that satisfies $\Phi(0)=\mathrm{Id}$
and $\Phi'=X\Phi$.

\begin{defin}
The \emph{Conley--Zehnder index} of the solution $z^{(1)}$ associated to the symplectic trivialization
$\Psi$, denoted by $\iCZ(z^{(1)},\Psi)$, is
the Conley--Zehnder of the path in $\Spl\big(T_{z(0)}\mathcal M,\varpi_{z(0)}\big)$ obtained by restriction
of the fundamental solution $\Phi$ to the interval $[0,T]$.
\end{defin}

\begin{rem}
It is known that, under suitable topological condition on the manifold $\mathcal M$ and
on the loop $z$, the Conley--Zehnder index $\iCZ(z^{(1)},\Psi)$ will not depend on the
choice of the trivialization $\Psi$. For instance, if $z:\mathds S^1\to \mathcal M$ is
homotopically trivial, i.e., if $z$ admits a continuous extension to the $2$-disk $\mathds D^2$,
then one can choose trivializations $\Psi$ of $z^*(T\mathcal M)$ that admit continuous
extensions to $\mathds D^2$. In this situation, if the first Chern class $c_1(\mathcal M)$
has vanishing integral on every $2$-sphere of $\mathcal M$, then the Conley--Zehnder
index of $z$ will not depend on $\Psi$. Namely, in this case any two trivializations
in the required class differ by a loop in the symplectic group which is homotopically
trivial which, by Corollary~\ref{thm:homtriloop}, does not alter the value of the Conley--Zehnder index.
\end{rem}

Using the results of Section~\ref{sec:prodformula}, we can estimate the Conley--Zehnder index of the iterated of
a periodic Hamiltonian solution as follows:
\begin{prop}\label{thm:iteratedhompower}
For all $N\ge1$, define $\Phi_{(N)}:[0,NT]\to\Spl\big(T_{z(0)}\mathcal M,\varpi_{z(0)}\big)$ to be the restriction to the interval $[0,NT]$
of the fundamental solution $\Phi$ of \eqref{eq:linhamequations}.
Then, $\Phi_{(N)}$ is fixed-endpoint homotopic to the $N$-th power of $\Phi_{(1)}$.
\end{prop}
\begin{proof}
Use the group law $\Phi(t+T)=\Phi(t)\cdot\Phi(T)$, and argue as in the proof of Lemma~\ref{thm:prodhomconc}.
\end{proof}

\begin{cor}
In the above notations, the following inequality holds:
\[\left\vert\iCZ\big(z^{(N)},\Psi\big)-N\cdot\iCZ\big(z^{(1)},\Psi\big)\right\vert\le\tfrac12n(N-1).
\]
In particular, $\left\vert\iCZ\big(z^{(N)},\Psi\big)\right\vert$ and $\left\vert\mu_{L_0}\big(z^{(N)},\Psi\big)\right\vert$ have sublinear growth in $N$;
moreover, if $\big\vert\iCZ\big(z^{(1)},\Psi\big)\big\vert>\tfrac12n$ (resp., if $\mu_{L_0}\big(z^{(1)},\Psi\big)>\tfrac72n$), then $\iCZ\big(z^{(N)},\Psi\big)$
(resp., $\mu_{L_0}\big(z^{(N)},\Psi\big)$) has linear growth
in $N$.
\end{cor}
\begin{proof}
Follows immediately from Corollary~\ref{thm:prodformula2} and Proposition~\ref{thm:iteratedhompower}.
\end{proof}

\end{section}

\section{The Weyl Representation of $\operatorname*{Mp}(2n,\mathbb{R)}$}
\label{sec:metaplectic}

\subsection{The metaplectic group}

We say that a quadratic form $W:\mathbb{R}^{n}\times\mathbb{R}^{n}%
\longrightarrow\mathbb{R}$ is \textquotedblleft
non-degenerate\textquotedblright\ if it can be written
\begin{equation}
W(x,x^{\prime})=\tfrac{1}{2}\left\langle Px,x\right\rangle -\left\langle
Kx,x^{\prime}\right\rangle +\tfrac{1}{2}\left\langle Qx^{\prime},x^{\prime
}\right\rangle \label{plq}%
\end{equation}
where $P$ and $Q$ are symmetric and $K$ invertible. The data of such a
quadratic form determine a symplectomorphism $\Phi_{W}\in\Spl(2n\R)$, whose matrix
in the canonical basis of $\R^{2n}$ is written in $n\times n$ blocks as:
\[
\begin{pmatrix}
K^{-1}Q\smallskip & K^{-1}\\
PK^{-1}Q-K^{T} & K^{-1}P
\end{pmatrix}
.\]
Set $L_0=\{0\}\times\R^n$; the symplectomorphism $\Phi_W$ can be characterized by the properties:
\begin{itemize}
\item $\Phi_W(L_{0})\cap L_{0}=\{0\}$
\item setting $z=(x,p),z'=(x',p')\in\R^{2n}$, then
$z=\Phi_{W}(z')$ if and only if  $p=\partial_{x}W(x,x^{\prime})$  and  $p'=-\partial_{x'}W(x,x^{\prime})$.
\end{itemize}

\begin{defin}
$\Phi_{W}$ is the \emph{free symplectic automorphism} determined by the \emph{generating
function $W$}.
\end{defin}

Let $\mathcal{S}(\mathbb{R}^{n})$ be the Schwartz space of rapidly decreasing
funcions on $\R^n$.  We associate to $\Phi_{W}$ the Fourier integral operator
$\widehat{\Phi}_{W,m}:\mathcal{S}(\R^{n})\longrightarrow
\mathcal{S}(\R^{n})$ defined by%
\[
\widehat{\Phi}_{W,m}f(x)=(2\pi i)^{-n/2}\Delta(W)\int e^{iW(x,x^{\prime}%
)}f(x^{\prime})d^{n}x^{\prime}
\]
where $\Delta(W)=i^{m}\sqrt{|\det K|}$ and $m$ corresponds to a choice of
$\arg\det K$ through
\[
m\pi=\arg(\det K)\mod4\pi
\]
(for each $K$ there are thus two choices of $m$ modulo 4). The operators
$\widehat{\Phi}_{W,m}$ extend by continuity to unitary operators on
$L^{2}(\R^{n});$ the inverse of $\widehat{\Phi}_{W,m}$ is
$\widehat{\Phi}_{W^{\ast},m^{\ast}}$ with $W^{\ast}(x,x^{\prime}%
)=-W(x^{\prime},x)$ and $m^{\ast}=n-m$. These operators thus generate a group
of unitary operators on $L^{2}(\R^{n})$, the \emph{metaplectic group}
$\operatorname*{Mp}(2n,\R)$, which is a double cover of
$\operatorname*{Sp}(2n,\R)$; the projection $\pi^{\operatorname*{Mp}%
}:\operatorname*{Mp}(2n,\R)\longrightarrow\Spl(2n,\R)$ is unambiguously determined
by the condition $\pi^{\operatorname*{Mp}}(\widehat{\Phi}_{W,m})=\Phi_{W}$. We have (see
\cite{mdglettmath}):

\begin{prop}
For every $\widehat{\Phi}\in\operatorname*{Mp}(2n,\R)$ there exist two
generating functions $W$ and $W^{\prime}$ and integers $m,m^{\prime}$ such
that $\widehat{\Phi}=\widehat{\Phi}_{W,m}\widehat{\Phi}_{W^{\prime},m^{\prime
}}$ and $\Phi_{W},\Phi_{W^{\prime}}\in\operatorname*{Sp}_{0}(2n,\R)$;
the condition $\Phi_{W}\in\operatorname*{Sp}_{0}(2n,\R)$ is equivalent
to $\det(P+Q-L-L^{T})\neq0$.
\end{prop}

The value modulo $4$ of $m+m^{\prime}-\mathrm n_-(P^{\prime}+Q)$
(recall the  $\mathrm n_-(R)$ denotes the index of the symmetric matrix $R$)
is independent of the choice of factorization, and thus only depends on
$\widehat{\Phi}$ (see \cite{AIF}):

\begin{defin}
The class modulo $4$ of $m+m^{\prime}-\mathrm n_-(P^{\prime}+Q)$ is
called the \emph{Maslov index} of $\widehat{\Phi}\in\operatorname*{Mp}(2n,\mathbb{R)}%
$; we denote it by $m(\widehat{\Phi})$. We call the function
$m:\operatorname*{Mp}(2n,\R)\longrightarrow\mathbb{Z}_{4}$ thus
defined "Maslov index" on $\operatorname*{Mp}(2n,\R)$.
\end{defin}

Let $\operatorname*{Inert}$ be the $2$-cocycle on $\Lambda$ defined by%
\begin{equation}
\operatorname*{Inert}(L,L^{\prime},L^{\prime\prime})=\tfrac{1}{2}%
(\tau(L,L^{\prime},L^{\prime\prime})+n+\partial\dim(L,L^{\prime}%
,L^{\prime\prime})\label{defi}%
\end{equation}
with $\partial\dim$ is the \v{C}ech coboundary of the $1$-cochain
$\dim(L,L^{\prime})=\dim(L\cap L^{\prime})$, that is
\[
\partial\dim(L,L^{\prime},L^{\prime\prime})=\dim(L\cap L^{\prime})-\dim(L\cap
L^{\prime\prime})+\dim(L^{\prime}\cap L^{\prime\prime}).
\]
We have (see \cite{AIF}):

\begin{prop}
We have $m(\widehat{\Phi}_{W,m})=\left[  m\right]  _{4}$ and the Maslov index
on $\operatorname*{Mp}(2n,\R)$ is related to the Leray index $\mu$ by
the formula%
\begin{equation}
m(\widehat{\Phi})=\left[  \tfrac{1}{2}(\mu(\Phi_{\infty}L_{0,\infty
},L_{0,\infty})+n+\dim(\Phi L_{0},L_{0}))\right]  _{4}\label{mu}%
\end{equation}
where $\Phi_{\infty}$ is any element of $\operatorname*{Sp}_{\infty
}(2n,\R)$ having projection $\widehat{\Phi}\in\operatorname*{Mp}%
(2n,\R)$ and $L_{0,\infty}$ is any element of $\Lambda_{\infty}$
covering $L_{0}=0\times\R^{n}$.
\end{prop}

It follows from formula (\ref{mu}) and the properties of the Leray index that%
\begin{equation}
m(\widehat{\Phi}\widehat{\Phi^{\prime}})=m(\widehat{\Phi})+m(\widehat
{\Phi^{\prime}})+\left[  \operatorname{Inert}(L_{0},\Phi L_{0},\Phi
\Phi^{\prime}L_{0})\right]  _{4}\text{.}\label{mfi}%
\end{equation}

\subsection{Weyl representation and Conley-Zehnder index}

Defining, as in \cite{WM}, the operator $R_{\nu}(\phi)$ associated to
$(\phi,\nu)\in\operatorname*{Sp}_{0}(2n,\R)\times\mathbb{Z}$ by the
Bochner integral%
\[
R_{\nu}(\phi)=\left(  \tfrac{1}{2\pi}\right)  ^{n}\frac{i^{\nu}}{\sqrt
{|\det(\phi-\mathrm{Id})|}}\int e^{\frac{i}{2}\left\langle M_{\phi}z,z\right\rangle
}T(z)d^{2n}z
\]
where $T(z)$ is the Heisenberg-Weyl operator, we have (see \cite{mdglettmath},
Prop. 6, \S 3.2 and Prop. 10, \S 3.3):

\begin{itemize}
\item Let $\phi_{W}$ be the free symplectic matrix generated by the quadratic
form (\ref{plq}). We have $\widehat{\Phi}_{W,m}=R_{\nu}(\phi_{W})$ if and only
$\nu=\nu(\widehat{\Phi}_{W,m})$ with%
\begin{equation}
\nu(\widehat{\Phi}_{W,m})\equiv m-\mathrm n_-(W_{xx})\mod4 \label{conca}%
\end{equation}
where $\mathrm n_-(W_{xx})$ is the index of inertia of the Hessian
matrix $W_{xx}$ of the function $x\longmapsto W(x,x)$;

\item Let $\widehat{\Phi}\in\operatorname*{Mp}(2n,\R)$ be such that
$\pi^{\operatorname*{Mp}}(\widehat{\Phi})\in\operatorname*{Sp}_{0}%
(2n,\R)$. If $\phi=\phi_{W}\phi_{W^{\prime}}$ and
\[
\widehat{\Phi}=R_{\nu(\widehat{\Phi}_{W,m})}(\phi_{W})R_{\nu(\widehat{\Phi
}_{W^{\prime},m^{\prime}})}(\phi_{W^{\prime}})
\]
then $\widehat{\Phi}=R_{\nu(\widehat{\Phi})}(\phi)$ with\footnote{
Recall that the Cayley transform used in this paper and that in \cite{mdglettmath} differ
by a sign.}
\begin{equation}
\nu(\widehat{\Phi})\equiv\nu(\widehat{\Phi}_{W,m})+\nu(\widehat{\Phi
}_{W^{\prime},m^{\prime}})-\tfrac{1}{2}\operatorname{sign}(M_{\phi_{W}%
}+M_{\phi_{W}^{\prime}})\mod4. \label{nunu}%
\end{equation}

\end{itemize}

These formulae suggest that there is a relation between the integer
$\nu(\widehat{\Phi})$ and the Conley-Zehnder index of some symplectic path
ending at $\phi=\pi_{\operatorname*{Mp}}(\widehat{\Phi})$. To study this
relation we will need the following two lemmas:

\begin{lem}
\label{lemsimple}(i) Let $(L,L^{\prime})\in(\Lambda(V,\omega))^{2}$. If $L\cap
L^{\prime\prime}=0$ then $\tau(L,L^{\prime},L^{\prime\prime})$ is the
signature of the quadratic form
\[
Q^{\prime}(z^{\prime})=\omega(\Pr\nolimits_{LL^{\prime\prime}}z^{\prime
},z^{\prime})=\omega(z^{\prime},\Pr\nolimits_{L^{\prime\prime}L}z^{\prime})
\]
on $L^{\prime}$, where $\Pr\nolimits_{LL^{\prime\prime}}$ is the projection
onto $L$ along $L^{\prime\prime}$ and $\Pr\nolimits_{L^{\prime\prime}L}%
=\mathrm{Id}-\Pr\nolimits_{LL^{\prime\prime}}$ is the projection on $L^{\prime\prime}$
along $L$. (ii) Let $(L,L^{\prime},L^{\prime\prime})$ be such that $L=L\cap
L^{\prime}+L\cap L^{\prime\prime}$. Then $\tau(L,L^{\prime},L^{\prime\prime
})=0$.
\end{lem}

\noindent(See \textit{e.g.} \cite{LV} for a proof).

\begin{lem}
\label{Serge}Let $\Phi:[0,1]\longrightarrow\operatorname*{Sp}(2n,\R%
)$\ be a continuous path such that $\Phi(0)=\mathrm{Id}$, $\Phi(1)=\phi\in
\operatorname*{Sp}_{0}(2n,\R)$. Then
\begin{equation}
\iCZ(\Phi)=-\tfrac{1}{2}\mu^{2}((\mathrm{Id}\oplus\phi)_{\infty}%
\Delta_{\infty},\Delta_{\infty}) \label{icz}%
\end{equation}
where $\Delta_{\infty}$ is any element of $\Lambda_{\infty}(\R%
^{4n},\omega^{2})$ with projection $\Delta=\{(z,z):z\in\R^{2n}\}$ and
$(\mathrm{Id}\oplus\phi)_{\infty}\in\operatorname{Sp}_{\infty}(\R^{4n}%
,\omega^{2})$ is the homotopy class in $\operatorname{Sp}(\R%
^{4n},\omega^{2})$ of the path $t\longmapsto\{(z,\Phi(t)z):z\in\R%
^{2n}\}$, $0\leq t\leq1$.
\end{lem}
\begin{proof} See \cite{PhD}. \end{proof}

It follows from the properties of the Leray
index (see \cite{JMPA}) that the right-hand side of (\ref{icz}) does not
depend on the choice of $\Delta_{\infty}$ covering $\Delta$.)

Let us now prove the main result of this section:

\begin{prop}
\label{Theodeux}Let $\Phi_{\infty}\in\operatorname*{Sp}_{\infty}%
(2n,\R)$ be the homotopy class of a continuous path $\Phi
:[0,1]\longrightarrow\operatorname*{Sp}(2n,\R)$\ such that $\Phi
(0)=\mathrm{Id}$, $\Phi(1)=\phi\in\operatorname*{Sp}_{0}(2n,\R)$. Let
$\widehat{\Phi}$ the image of $\Phi_{\infty}$ in $\operatorname*{Mp}%
(2n,\R)$ by the covering mapping $\operatorname*{Sp}_{\infty
}(2n,\R)\longrightarrow\operatorname*{Mp}(2n,\R)$. We have%
\begin{equation}
\nu(\widehat{\Phi})\equiv-\iCZ(\Phi)\text{ \ }%
\operatorname{mod}4. \label{apo}%
\end{equation}

\end{prop}

\begin{proof}
Since $\widehat{\Phi}$ can be written as a product $\widehat{\Phi}%
_{W,m}\widehat{\Phi}_{W^{\prime},m^{\prime}}$, formula (\ref{nunu}) and the
product formula~\eqref{eq:firstprodform} in Corollary~\ref{thm:prodformula}
reduce the proof to the case $\widehat{\Phi
}=\widehat{\Phi}_{W,m}$. In view of Lemma \ref{Serge} and (\ref{conca}) it
is sufficient to show that
\begin{equation}
m-\mathrm n_-(W_{xx})\equiv\tfrac{1}{2}\mu^{2}((\mathrm{Id}\oplus\phi
_{W})_{\infty}\Delta_{\infty},\Delta_{\infty})\text{\ }\operatorname{mod}4.
\label{apa}%
\end{equation}
We will divide the proof of (\ref{apa}) in three steps. We denote as before by
$\omega^{2}$ the symplectic form $\omega\oplus(-\omega)$ on $\R%
^{4n}=\R^{2n}\oplus\R^{2n}$; the corresponding Wall--Kashiwara
and Leray indexes are $\tau^{2}$ and $\mu^{2}$.

\noindent
\textbf{(I)} Let $L^{(2)}%
\in\Lambda(\R^{4n},\omega^{2})$. Let $L_{\infty}^{(2)}\in
\Lambda_{\infty}(\R^{4n},\omega^{2})$ cover $L^{(2)}\in\Lambda
(\R^{4n},\omega^{2})$. Using the property $\partial\mu^{2}=\pi^{\ast
}\tau^{2}$ we get after a few calculations%
\begin{multline}
\mu^{2}((\mathrm{Id}\oplus\phi_{W})_{\infty}\Delta_{\infty},\Delta_{\infty})=(\mu
^{2}((\mathrm{Id}\oplus\phi_{W})_{\infty}L_{\infty}^{(2)},L_{\infty}^{(2)}%
)\label{star}\\
+\tau^{2}((\mathrm{Id}\oplus\phi_{W})\Delta,\Delta,L^{(2)})-\tau^{2}((\mathrm{Id}\oplus\phi
_{W})\Delta,(\mathrm{Id}\oplus\phi_{W})L^{(2)},L^{(2)}))\text{.}\nonumber
\end{multline}
Choosing in particular $L^{(2)}=L_{0}^{(2)}=L_{0}\oplus L_{0}$ (recall: $L_0=\{0\}\times\R^n$)
we get
\begin{align*}
\mu^{2}((\mathrm{Id}\oplus\phi_{W})_{\infty}L_{0,\infty}^{(2)},L_{0,\infty}^{(2)})  &
=\mu^{2}((\mathrm{Id}\oplus\phi_{W})_{\infty}(L_{0}\oplus L_{0})_{\infty},(L_{0}\oplus
L_{0})_{\infty})\\
&  =\mu(L_{0,\infty},L_{0,\infty})-\mu(L_{0,\infty},\Phi_{W,\infty}%
L_{P,\infty})\\
&  =\mu(\Phi_{W,\infty}L_{0,\infty},L_{0,\infty})
\end{align*}
so that there remains to prove that
\[
\tau^{2}((\mathrm{Id}\oplus\phi_{W})\Delta,\Delta,L_{0}^{(2)})-\tau^{2}((\mathrm{Id}\oplus\phi
_{W})\Delta,(\mathrm{Id}\oplus\phi_{W})L_{0}^{(2)},L_{0}^{(2)})=-2\operatorname*{sign}%
W_{xx}\text{.}%
\]
\textbf{(II)} We are going to show that $\tau^{2}((\mathrm{Id}\oplus\phi_{W}%
)\Delta,(\mathrm{Id}\oplus\phi_{W})L_{0}^{(2)},L_{0}^{(2)})=0$; in view of the
symplectic invariance and the antisymmetry of $\tau^{2}$ this is equivalent to%
\begin{equation}
\tau^{2}(L_{0}^{(2)},\Delta,L_{0}^{(2)},(\mathrm{Id}\oplus\phi_{W})^{-1}L_{0}%
^{(2)})=0\text{.} \label{cucu}%
\end{equation}
We have
\[
\Delta\cap L_{0}^{(2)}=\{(0,p;0,p):p\in\R^{n}\}
\]
and $(\mathrm{Id}\oplus\phi)^{-1}L_{0}^{(2)}\cap L_{0}^{(2)}$ consists of all
$(0,p^{\prime},\phi^{-1}(0,p^{\prime\prime}))$ with $\phi^{-1}(0,p^{\prime
\prime})=(0,p^{\prime})$; since $\phi_{W}$ (and hence also $\phi_{W}^{-1}$) is
free we must have $p^{\prime}=p^{\prime\prime}=0$ so that%
\[
(\mathrm{Id}\oplus\phi_{W})^{-1}L_{0}^{(2)}\cap L_{0}^{(2)}=\{(0,p;0,0):p\in
\R^{n}\}\text{.}%
\]
It follows that we have%
\[
L_{0}=\Delta\cap L_{0}^{(2)}+(\mathrm{Id}\oplus\phi_{W})^{-1}L_{0}^{(2)}\cap L_{0}^{(2)}%
\]
hence (\ref{cucu}) in view of property \textit{(ii)} in Lemma \ref{lemsimple}.
\smallskip

\noindent
\textbf{(III) }Let us finally show that.%
\[
\tau^{2}((\mathrm{Id}\oplus\phi_{W})\Delta,\Delta,L_{0}^{(2)})=-2\operatorname*{sign}%
W_{xx}\text{;}%
\]
this will complete the proof. The condition $\det(\phi_{W}-\mathrm{Id})\neq0$ is
equivalent to $(\mathrm{Id}\oplus\phi_{W})\Delta\cap\Delta=0$ hence, using property
\textit{(i)} in Lemma \ref{lemsimple}:%
\[
\tau^{2}((\mathrm{Id}\oplus\phi_{W})\Delta,\Delta,L_{0}^{(2)})=-\tau^{2}(\mathrm{Id}\oplus\phi
_{W})\Delta,L_{0}^{(2)},\Delta)
\]
is the signature of the quadratic form $Q$ on $L_{0}$ defined by
\[
Q(0,p,0,p^{\prime})=-\omega^{2}(\Pr\nolimits_{(\mathrm{Id}\oplus\phi_{W})\Delta,\Delta
}(0,p,0,p^{\prime});0,p,0,p^{\prime})
\]
where
\[
\Pr\nolimits_{(\mathrm{Id}\oplus\phi_{W})\Delta,\Delta}=%
\begin{bmatrix}
(\phi_{W}-\mathrm{Id})^{-1} & -(\phi_{W}-\mathrm{Id})^{-1}\smallskip\\
\phi_{W}(\phi_{W}-\mathrm{Id})^{-1} & -\phi_{W}(\phi_{W}-\mathrm{Id})^{-1}%
\end{bmatrix}
\]
is the projection on $(\mathrm{Id}\oplus\phi_{W})\Delta$ along $\Delta$ in
$\R^{2n}\oplus\R^{2n}$. It follows that the quadratic form $Q$
is given by%
\[
Q(0,p,0,p^{\prime})=-\omega^{2}((\mathrm{Id}-\phi_{W})^{-1}(0,p^{\prime\prime}),\phi
_{W}(\mathrm{Id}-\phi_{W})^{-1}(0,p^{\prime\prime});0,p,0,p^{\prime})
\]
where we have set $p^{\prime\prime}=p-p^{\prime}$; by definition of
$\omega^{2}$ this is%
\[
Q(0,p,0,p^{\prime})=-\omega((\mathrm{Id}-\phi_{W})^{-1}(0,p^{\prime\prime}%
),(0,p))+\omega(\phi_{W}(\mathrm{Id}-\phi)^{-1}(0,p^{\prime\prime}),(0,p^{\prime
}))\text{. }%
\]
We have%
\[
(\mathrm{Id}-\phi_{W})^{-1}=JM_{\phi_{W}}+\tfrac{1}{2}\mathrm{Id}\text{ \ , \ }\phi_{W}(\mathrm{Id}-\phi
_{W})^{-1}=JM_{\phi_{W}}-\tfrac{1}{2}\mathrm{Id}
\]
and hence%
\begin{align*}
Q(0,p,0,p^{\prime})  &  =-\omega((JM_{\phi_{W}}+\tfrac{1}{2}\mathrm{Id})(0,p^{\prime
\prime}),(0,p))+\omega((JM_{\phi_{W}}-\tfrac{1}{2}\mathrm{Id})(0,p^{\prime\prime
}),(0,p^{\prime}))\\
&  =-\omega(JM_{\phi_{W}}(0,p^{\prime\prime}),(0,p))+\omega(JM^{\phi_{W}%
}(0,p^{\prime\prime}),(0,p^{\prime}))\\
&  =\omega(JM_{\phi_{W}}(0,p^{\prime\prime}),(0,p^{\prime\prime}))\\
&  =-\left\langle M_{\phi_{W}}(0,p^{\prime\prime}),(0,p^{\prime\prime
})\right\rangle \text{.}%
\end{align*}
Let us calculate explicitly $M_{\phi_{W}}$. Writing $\phi_{W}=%
\begin{bmatrix}
A & B\\
C & D
\end{bmatrix}
$ we have%
\[
\phi_{W}-\mathrm{Id}=%
\begin{bmatrix}
0 & B\smallskip\\
\mathrm{Id} & D-\mathrm{Id}
\end{bmatrix}%
\begin{bmatrix}
C-(D-\mathrm{Id})B^{-1}(A-\mathrm{Id}) & 0\smallskip\\
B^{-1}(A-\mathrm{Id}) & \mathrm{Id}
\end{bmatrix}
\]
that is%
\begin{equation}
\phi_{W}-\mathrm{Id}=%
\begin{bmatrix}
0 & B\\
\mathrm{Id} & D-\mathrm{Id}
\end{bmatrix}%
\begin{bmatrix}
W_{xx} & 0\\
B^{-1}(A-\mathrm{Id}) & \mathrm{Id}
\end{bmatrix}
\label{smoinsi}%
\end{equation}
where we have used the identity
\[
C-(D-\mathrm{Id})B^{-1}(A-\mathrm{Id}))=B^{-1}A+DB^{-1}-B^{-1}-(B^{T})^{-1}%
\]
which follows from the relation $C-DB^{-1}A=-(B^{T})^{-1}$ due to the fact
that $\phi_{W}$ is symplectic. We thus have, writing $W_{xx}^{-1}%
=(W_{xx})^{-1}$,
\begin{align*}
(\phi_{W}-\mathrm{Id})^{-1}  &  =%
\begin{bmatrix}
W_{xx}^{-1} & 0\smallskip\\
B^{-1}(\mathrm{Id}-A)W_{xx}^{-1} & \mathrm{Id}
\end{bmatrix}%
\begin{bmatrix}
(\mathrm{Id}-D)B^{-1} & \mathrm{Id}\smallskip\\
B^{-1} & 0
\end{bmatrix}
\medskip\\
&  =%
\begin{bmatrix}
W_{xx}^{-1}(\mathrm{Id}-D)B^{-1} & W_{xx}^{-1}\smallskip\\
B^{-1}(\mathrm{Id}-A)W_{xx}^{-1}(\mathrm{Id}-D)B^{-1}+B^{-1} & B^{-1}(\mathrm{Id}-A)W_{xx}^{-1}%
\end{bmatrix}
\end{align*}
and hence%
\[
M_{\phi_{W}}=%
\begin{bmatrix}
B^{-1}(\mathrm{Id}-A)W_{xx}^{-1}(\mathrm{Id}-D)B^{-1}+B^{-1} & \frac{1}{2}\mathrm{Id}+B^{-1}(\mathrm{Id}-A)W_{xx}%
^{-1}\smallskip\\
-\frac{1}{2}\mathrm{Id}-W_{xx}^{-1}(\mathrm{Id}-D)B^{-1} & -W_{xx}^{-1}%
\end{bmatrix}
\]
so that we have%
\[
Q(0,p,0,p^{\prime})=\left\langle W_{xx}^{-1}p^{\prime\prime},p^{\prime\prime
}\right\rangle =\left\langle W_{xx}^{-1}(p-p^{\prime}),(p-p^{\prime
})\right\rangle \text{.}%
\]
It follows that the matrix of the quadratic form $Q$ is%
\[
2%
\begin{bmatrix}
W_{xx}^{-1} & -W_{xx}^{-1}\smallskip\\
-W_{xx}^{-1} & W_{xx}^{-1}%
\end{bmatrix}
\]
and this matrix has signature $2\operatorname*{sign}(W_{xx})^{-1}%
=2\operatorname*{sign}W_{xx}$. This concludes the proof since $\det(\phi
_{W}-\mathrm{Id})=(-1)^{n}\det B\det W_{xx}$ and $2m-m=\mu(\phi_{W,\infty}L_{0,\infty
},L_{0,\infty})$.
\end{proof}

\end{document}